\documentclass[a4paper,10pt]{article}

\usepackage{color}
\usepackage{amsmath}
\usepackage{amsfonts}
\usepackage{amssymb}
\usepackage{graphicx}
\usepackage{epsfig}
\newtheorem{montheo}{Theorem}

\newtheorem{prop}{Proposition}
\newtheorem {lem}{Lemma}
\newtheorem {rem}{\textit{Remark}}

\title{Quadratic centers defining Elliptic Surfaces}
\author{S\'ebastien GAUTIER}

\begin{document}

\maketitle

\begin{abstract}
Let $X$ be a quadratic vector field with a center whose generic orbits are algebraic curves of genus one. To each $X$ we associate an elliptic surface (a smooth complex compact surface which is a genus one fibration). We give the list of all such vector fields and determine the corresponding elliptic surfaces.

\end{abstract}

\section{Introduction }
\indent The second part of the $16$th Hilbert problem asks for an upper bound  to the number of limit cycles of a plane
 polynomial vector field of degree less or equal to $n$.
 Even in the case of quadratic systems ($n=2$) the problem remains open. An infinitesimal version of the  $16$th
 Hilbert problem can be formulated as follows:

 {\em Find an upper bound $Z (f,n) $ to the number of limit cycles of a polynomial vector field of degree $n$,
 close to a polynomial vector field with a first integral $f$}.

 The associated foliation on the plane is defined by
\begin{equation}\label{perturbation}
 R^{-1}df + \varepsilon \omega =0
\end{equation}
 where $R^{-1}df= Pdx +Qdy$ is a given polynomial one-form, $\deg P, \deg Q \leq n$, $R^{-1}$ is an integrating factor,
 and $\omega$ is a polynomial one-form of degree $n$ with coefficients depending analytically on the small parameter
 $\varepsilon$.

 A progress in solving the infinitesimal 16th Hilbert problem is achieved mainly in  the case when  $ f $ is a
 polynomial of degree three,
 or $F=y^2+P(x)$ where $P$ is a polynomial of degree four see \cite{Ilyashenko, Petrov, hamiltonian}.
A key point is that  the generic leaves  $\Gamma_c = \{f=c\}\subset \mathbb{C}^2$ of the
polynomial foliation
 $R^{-1}df= 0$ are
 elliptic curves. We expect that the perturbations of more general polynomial foliations with elliptic leaves
 (which we call "elliptic foliations") can be studied along the same lines. This leads naturally to the following (open)
 problem.

\emph{For a given $n>1$ determine, up to an affine equivalence, the elliptic polynomial foliations $Pdx+Qdy=0$,
$\deg P, \deg Q \leq n$.}

The present paper adresses the above problem in the quadratic case, $n=2$. In view of  applications to the
$16$th Hilbert problem most important is the case when the non-perturbed foliation is real and possesses a
center. Such foliations are well-known since Dulac (1908)  and Kapteyn
(1912). Moreover, when the leaves of the foliation (the orbits of the quadratic vector
field) are algebraic curves, there is a (rational) first integral
$f$ \cite{pfaff}. Reminding the classification of quadratic vector fields with a
center, a rational first integral of the foliation induced is thus of
four different kinds:
\begin{eqnarray}
f=P_3(x,y) \ \textrm{with} \ P_3 \in \mathbb{R}_3[x,y] \ \mbox{\quad (Hamiltonian case)}
\end{eqnarray}
\begin{eqnarray}
\label{rev} f=x^{\lambda}(y^{2}+\mathit{P}_{2}(x)) \ \textrm{with} \ 
\left \lbrace 
\begin{array}{l}
\lambda \in \mathbb{Q} \\
 P_2 \in \mathbb{R}_2[x,y]\\
\end{array}
\right. \mbox{\quad  (reversible case)}
\end{eqnarray}
\begin{eqnarray}
\label{LV} f=x^{\lambda}y^{\mu}(ax+by+c)\ \textrm{with} \ \left \lbrace 
\begin{array}{l}
\lambda, \mu \in \mathbb{Q} \\
 P_2 \in \mathbb{R}_2[x,y]\\
$a$, $b$, $c$ \ \textrm{real numbers}\\
\end{array}
\right. \mbox{(Lotka-Volterra case)}
\end{eqnarray}
\begin{eqnarray}
\label{c4}
f=P_2(x,y)^{-3}P_3(x,y)^2 \ \textrm{with} \
\left \lbrace 
\begin{array}{l}
 P_3 \in \mathbb{R}_3[x,y] \\
 P_2 \in \mathbb{R}_2[x,y]\\
\end{array}
\right.  \ \ 
\mbox{(codimension $4$ case)}
\end{eqnarray}

In section \ref{section2} we give the classification, up to an affine equivalence, of all elliptic foliations with
a first integral of the form (\ref{rev}) or (\ref{LV}). The
Hamiltonian cases obviously induce an elliptic foliation and have
already been studied. Remarks concerning the codimension $4$ case can
be found in \cite{ggi, thesis}. In our classification, the base
field is supposed to be $\mathbb{C}$, so all parameters $a,b,c,\lambda,...$  are complex. We get a finite list
of such foliations with a center as well several infinite series of degenerated foliations which can not have a
center (when the base field is $\mathbb{R}$). Most of these  elliptic foliations  were not previously studied in
the context of the 16th Hilbert problem (but see \cite{chen,yu,
  iliev}).\\
\\
\indent The section \ref{sec3} deals with the topology of the singular
surface induced. An elliptic foliation in $\mathbb{C}^2$ (or more generally a foliation with an algebraic first integral $f$)
gives rise canonically to an elliptic surface as follows. Suppose that $f$ is chosen in such a way that the
generic fiber of the map $f:\mathbb{C}^2\rightarrow \mathbb{C}$ is an irreducible algebraic curve. The
induced
 rational map $f: \mathbb{P}^2\dasharrow \mathbb{P}^1$ have a finite number of points of indetermination.
After a finite number of blow-ups of $\mathbb{P}^2$ at these points we get (by the Hironaka's desingularisation
theorem ) an induced analytic map
$$\mathbb{P}^2 \subset K \stackrel{f}{\rightarrow} \mathbb{P}^1$$
where $K$ is a smooth complex surface.
 We may further suppose that $K$ is minimal in the sense that the
fibers do not contain exceptional curves of first kind. The pair $(K,f)$ is then the elliptic surface associated
to the elliptic foliation $R^{-1}df= 0$. It is unique up to a fiber
preserving isomorphism. In this last section, we compute the singular
fibers of the elliptic surfaces obtained. The singular fibers of an
elliptic surface are classified by Kodaira \cite{singularbis}. Such 
computations are $2$-folds. First of all, it permits to identifiate isomorphic
elliptic surfaces of non affine equivalent foliations, wich on its own
is of interest. But the most important is that it immediately gives
the local monodromy of the singular fibers. Here, the number of
singular fibers (except in the Hamiltonian case) do not exceed $4$ so
the local monodromy of the singular fibers gives a good description
of  the (global) monodromy group of the associated Picard-Fuchs equation (or
equivalently, the homological invariant of the surface \cite{singularbis}), which on its turn is necessary when
studying zeros of Abelian integrals (or limit cycles of the perturbed foliation (\ref{perturbation})), see
\cite{hamiltonian,Petrov,Ilyashenko,ggi, thesis} for details. 
\section{Quadratic centers which define elliptic foliations  }
\label{section2}
Let $\mathcal{F} = \mathcal{F}(\omega)$ be a foliation on the plane $\mathbb{C}^2$ defined by a differential
form $\omega = Pdx +Q dy$. We say that $\mathcal{F}(df)$ is \emph{elliptic} provided that its generic leaves are
elliptic curves.
 As stated in
the Introduction, in the present paper we suppose that $\mathcal{F}$ has, eventually after an affine change of
the variables in $\mathbb{C}^2$, a first integral of the form (\ref{rev}) or (\ref{LV}). Such a foliation will
be called \emph{reversible} (having an integral of the form (\ref{rev}) but not (\ref{LV})), of \emph{
Lotka-Voltera type} (having an integral of the form (\ref{LV}) but not (\ref{rev})), or  of reversible
Lotka-Voltera type.

\subsection{The  reversible case}
\label{section21}
An elliptic foliation of Lotka-Voltera type has three invariant lines. From this we deduce that a reversible
Lotka-Voltera foliation has always a first integral $f=x^{\lambda}(y^{2}+\mathit{P}_{2}(x))$ where $P_2$ is a
polynomial of degree at most two, and the bi-variate polynomial $y^{2}+\mathit{P}_{2}(x)$ is irreducible. In
this section we prove the following:
\begin{montheo}
\label{th1} The reversible foliation $\mathcal{F}(df)$ is elliptic if and only if, after an affine change of the
variables, it has a first integral of the form:
$$\begin{array}{ll}  {\rm (rv1)}\; f=x^{-3}(
y^2+ax^2+bx+c) & {\rm (rv2)}\; f=x( y^2+cx^2+bx+a)\\
{\rm (rv3)}\; f=x^{-3/2}( y^2+ax^2+bx+c) & {\rm (rv4)}\; f=x^{-1/2}( y^2+cx^2+bx+a)\\ {\rm
(rv5)}\; f=x^{-4}(
y^2+ax^2+bx+c) & {\rm (rv6)}\; f=x^2(y^2+cx^2+bx+a)\\
{\rm (rv7)}\; f=x^{-4/3}(y^2+bx+c) & {\rm (rv8)}\; f=x^{-2/3}( y^2+cx^2+bx)\\ {\rm (rv9)}\;
f=x^{-4/3}( y^2+ax^2+bx) & {\rm (rv10)}\; f=x^{-2/3}( y^2+bx+a)\\ {\rm (rv11)}\; f=x^{-5/3}(
y^2+ax^2+bx) & {\rm (rv12)}\; f=x^{-1/3}( y^2+bx+a)\\ {\rm (rv13)}\; f=x^{-5/4}( y^2+ax^2+bx)
& {\rm (rv14)}\; f=x^{-3/4}( y^2+bx+a)\\ {\rm (rv15)}\; f=x^{-7/4}(
y^2+ax^2+bx) & {\rm (rv16)}\; f=x^{-1/4}( y^2+bx+a)\\
{\rm (rv17)}\; f=x^{-5/2}( y^2+ax^2+bx) & {\rm (rv18)}\; f=x^{1/2}( y^2+bx+a).
\end{array}$$
or
$$
(i) f=x^{-1 + \frac{2}{k}}(y^{2}+x), k\in\mathbb{Z}^*\setminus 2\mathbb{Z}, $$
$$
(ii) f=x^{-1 + \frac{3}{k}}(y^{2}+x), k\in
\mathbb{Z}^* \setminus 3\mathbb{Z}.$$
\end{montheo}
\begin{rem}
\label{Remarque1}
We shall assume moreover that $c \neq 0$ for $\rm (rv3), \rm(rv4).$
\end{rem}
{\bf Proof.} Let  $\Gamma_t$ be the set of $(x,y)\in \mathbb{C}^2$ such that for some determination of the
multi-valued function $x^{\lambda}$ holds $f(x,y)=t $. If the connected components of $\Gamma_t$ for all $t$ are
algebraic curves, then $\lambda \in \mathbb{Q}$ and we put $\lambda=\frac{p}{q}$ ($p \in \mathbb{Z}$, $q \in
\mathbb{N}^{*}$ and $gcd(p,q)=1$) with $\mathit{P}_{2}(x)=ax^{2}+bx+c \in \mathbb{C}_{2}[x].$

As the foliation is reversible  we may suppose that $y^2+ax^{2}+bx+c$ is irreducible, or simply $b^2-4ac\neq 0$.
We shall  suppose first that $a\neq 0, c\neq 0$, that is to say the quadric $\{y^2+ax^{2}+bx+c=0\}$ is not
tangent to the line at infinity in $\mathbb{P}^2$ and to the line $\{x=0\}$.
\subsubsection{The case $a\neq 0, c\neq 0,
b^2-4ac\neq 0$.}
 After a scaling of $t$ and an affine transformation we may suppose that
\begin{eqnarray}
f= x^{\lambda}(y^{2}+x^{2}+bx+c).
\end{eqnarray}
By abuse of notation we put
$$
\Gamma_t=\{x^{p/q}(y^{2}+x^{2}+bx+c)=t \}
$$
and in a similar way we define
\begin{eqnarray}
\label{forme1} \tilde{\Gamma}_t=\{ X^{p}(Y^{2}+X^{2q}+bX^{q}+c)=t\}.
\end{eqnarray}
\begin{lem}
\label{lem1}
The map $\varphi : \mathbb{C}^2 \rightarrow \mathbb{C}^2 : (X,Y) \rightarrow (x,y)=(X^q,Y)$ induces an
isomorphism of $\Gamma_t$ and $\tilde{\Gamma}_t$.
\end{lem}
Indeed, it is straightforward to check that $\varphi : \tilde{\Gamma}_t \rightarrow \Gamma_t$ is a bijection and
therefore is a bi-holomorphic map.

To compute the genus of $\tilde{\Gamma}_t$ or $\Gamma_t$ we distinguish two cases:

\begin{enumerate}
 \item The case when $p<0.$
\\
\\
We obtain the hyper-elliptic curve $\{ y^{2}=-x^{2q}-bx^{q}+tx^{-p}-c\}$. The roots of the polynomial
$-x^{2q}-bx^{q}+tx^{-p}-c$ are different and non zeros since $t$ is generic. Consequently, its genus is one if and only if the
degree of the polynomial is $3$ or $4$ and thus we get:
\begin{enumerate}
\item
$  f = x^{-3}(y^{2}+x^{2}+bx+c) $
\item
$  f= x^{-4}(y^{2}+x^{2}+bx+c) $
\item
$  f = x^{-\frac{1}{2}}(y^{2}+x^{2}+bx+c) $
\item
$  f = x^{-\frac{3}{2}}(y^{2}+x^{2}+bx+c) $
\end{enumerate}

 \item Suppose now $p\geq0.$
\\
\\
We easily have: $y^{2}x^{p}=t-x^{2q+p}-bx^{q+p}-cx^{p}.$
Thus after a birational transformation, we obtain:
$$y^{2}=x^{p}(t-x^{2q+p}-bx^{q+p}-cx^{p}).$$
Since $t$ is generic, all the roots of $t-x^{2q+p}-bx^{q+p}-cx^{p}$ are different and do not vanish.
\begin{itemize}
 \item If $p$ is even we have:
$$(\dfrac{y}{x^{\frac{p}{2}}})^{2}=t-x^{2q+p}-bx^{q+p}-cx^{p}.$$
Consequently, it is elliptic when $2q+p$ either equal $3$ or $4$. This gives the following curves:
\\
\begin{enumerate}
\item
$x(y^{2}+x^{2}+bx+c) = t$
\item
$x^{2}(y^{2}+x^{2}+bx+c) = t.$
\end{enumerate}
 \item If $p$ is odd, then we have:
$$(\dfrac{y}{x^{\frac{p-1}{2}}})^{2}=x(t-x^{2q+p}-bx^{q+p}-cx^{p}).$$
Since all the roots of $t-x^{2q+p}-bx^{q+p}-cx^{p}$ are different and non zeros, the curve is elliptic if and only
$2q+p$ either equals $2$ or $3$, which gives the solution $\rm(b)$ above.
\end{itemize}
\end{enumerate}
This we have obtained the cases (rv1)-(rv6) in Theorem \ref{th1}.
\subsubsection{The case $a\neq 0, c= 0,
b^2-4ac\neq 0$.} This means that the quadric $\{y^2+ax^{2}+bx=0\}$ is tangent to the line $\{x=0\}$ and is
transversal to the line at infinity, see Figure \ref{example2}. After an affine transformation and scaling of
$t$ we get  $P_2(x)= ax^{2}+bx$ with $a  \neq 0$. Therefore we need to compute  the genus of
$\{x^{p}(y^{2}+x^{2q}+bx^{q})=t\}$ for generic $t$.\\
\\
If $p \geq 0$ the same computations as  case $a\neq 0, c\neq0$ give the same  solutions of the problem.
\\
Let $p<0$ and suppose that $p=2a$ is even.
We have $(x^{a}y)^{2}=-x^{2q+p}-bx^{q+p}+t$, so if $-p\leq q$, it has genus one if and only $2q+p=3$ or $4$,
hence $q\leq 4$ and $(p,q)=(-2,3)$.\\
If $-p \geq 2q$ the curve above is birational  to the curve $y^{2}=-x^{-p-q}-bx^{-p-2q}+t$ and so
it has genus one if $-p-q=3$ or $4$ which leads to $(p,q)=(-4,1).$\\
If $q<-p<2q$, because $p$ is even and $q$ is odd, it is equivalent to calculate the genus of
$\{y^{2}=x(tx^{-(q+p)}-x^{q}-b)\}$. Therefore $(p,q)=(-4,3)$.\\
\\
Now suppose that $p=2a+1$ is odd. The curve is birationally equivalent  to $\{y^{2}=-x(x^{2q+p}+bx^{q+p}-t)\}$.
As above we get: $(p,q)=(-1,2)$, $(-3,1)$,$(-5,2)$, $(-5,3)$, $(-5,4)$, $(-7,4). $

To resume, we proved
\begin{prop}
The  foliation $\mathcal{F}(x^{\lambda}(y^{2}+x^{2}+bx))$ is elliptic if and only it has under affine transformation a first integral of the kind : \\

\begin{tabular}{l l l l l l }

$f=x(y^{2}+x^{2}+bx) $ &
$f=x^{2}(y^{2}+x^{2}+bx) $ \\
\\
$f= x^{-3}(y^{2}+x^{2}+bx) $ &
$f= x^{-4}(y^{2}+x^{2}+bx) $\\
\\
$f= x^{-\frac{1}{2}}(y^{2}+x^{2}+bx) $ &
$f= x^{-\frac{5}{2}}(y^{2}+x^{2}+bx) $\\
\\
$f= x^{-\frac{2}{3}}(y^{2}+x^{2}+bx) $ &
$f= x^{-\frac{4}{3}}(y^{2}+x^{2}+bx) $ \\
\\
$f= x^{-\frac{5}{4}}(y^{2}+x^{2}+bx) $ &
$f= x^{-\frac{7}{4}}(y^{2}+x^{2}+bx) $ \\
\\
$f= x^{-\frac{5}{3}}(y^{2}+x^{2}+bx) $ &

\end{tabular}

\end{prop}
Here we get $(\rm rv1)-(\rm rv6)$ except $(\rm rv3)$ according to Remark \ref{Remarque1}, $\rm(rv8)$ and the end of the left column of Theorem \ref{th1}.

\subsubsection{ $a=0, c\neq 0$, $b^2-4ac\neq 0$.}
\label{213}
This means that the quadric $\{y^2+bx+c=0\}$ is tangent to the line at infinity and is transversal to the
line $\{x=0\}$. The birational change of variables $x\rightarrow 1/x$, $y\rightarrow y/x$ shows that this is
equivalent to the case $a\neq 0, c=0$ and we get:
\begin{prop}
The  foliation $\mathcal{F}(x^{\lambda}(y^{2}+x+c))$ is elliptic if and only it has under affine transformation a first integral of the kind : \\
\\
\begin{tabular}{l l l l l l }

$f=x(y^{2}+x+c) $ &
$f=x^{2}(y^{2}+x+c) $ \\
\\
$ f=x^{-3}(y^{2}+x+c) $ &
$ f=x^{-4}(y^{2}+x+c) $\\
\\
$ f=x^{-\frac{3}{2}}(y^{2}+x+c) $ &
$ f=x^{\frac{1}{2}}(y^{2}+x+c) $\\
\\
$f= x^{-\frac{2}{3}}(y^{2}+x+c) $ &
$f= x^{-\frac{4}{3}}(y^{2}+x+c) $ \\
\\
$f= x^{-\frac{1}{4}}(y^{2}+x+c) $ &
$f=x^{-\frac{3}{4}}(y^{2}+x+c) $ \\
\\
$f=x^{-\frac{1}{3}}(y^{2}+x+c) $ &

\end{tabular}

\end{prop}
Here we get $(\rm rv1)-(\rm rv6)$ except $(\rm rv4)$ according to Remark \ref{Remarque1}, $\rm(rv7)$ and the end of the right column of Theorem \ref{th1}.
\subsubsection{$a=c=0, b\neq 0$}
This means that the quadric $\{y^2+ax^{2}+bx+c=0\}$ is tangent to the line at infinity and  to the line
$\{x=0\}$.
 Up to affine change of re-scalings we may suppose  $f=x^{p/q}(y^{2}+x)$. \\
\\
If $p$  is even the curve $\tilde{\Gamma}_t $ is birational to $y^{2}=-x^{q+p}+t.$\\If $p $ is odd the curve $\tilde{\Gamma}_t$ is birational to $y^{2}=-x(x^{q+p}-t).$\\
For $p\geq 0$ this curve is elliptic if and only: $(p,q)=(1,1)$, $(2,1)$ and $(1,2)$.
\\
Now, if $p\leq 0$ and $q+p\geq0$ the conditions are $q+p=2 \ ,\ 3$ or $4$ with $p$ even.\\
The case $q+p\leq0$ gives similarly$-q-p=2 \ , \ 3$ or $4$ with $p$ even.\\
Notice that we must have q prime with the integers $2$ or $3$ or $4$ when considering all cases. This gives the following:
\begin{prop}
The  foliation $\mathcal{F}(x^{\lambda}(y^{2}+bx))$ with $b \neq 0$ is elliptic  if and only if it has a first integral of the kind : \\
\\
\begin{tabular}{lc l c l }

$f=x^{-1+ \frac{2}{k}}(y^{2}+x)$, $k \in \mathbb{Z}^{*}\setminus 2 \mathbb{Z}$& $f=x^{-1+ \frac{3}{k}}(y^{2}+x)$, $k \in \mathbb{Z}^{*}\setminus 3 \mathbb{Z}$
\\

\end{tabular}

\end{prop}
Finally, Theorem \ref{th1} is proved. 

\begin{flushright}
$\square$
\end{flushright}

\subsection{The  Lotka-Volterra case}
\label{section22}
\begin{montheo} 
\label{th2} The Lotka-Volterra foliation $\mathcal{F}(df)$ is elliptic if
and only if, after an affine change of the variables, it has a first integral of the form:\\

$$\begin{array}{lll}
{\rm (lv1)}\; f=x^{2}y^{3}(1-x-y)  &
{\rm (lv2)}\; f=x^{-6}y^{2}(1-x-y)  &
{\rm (lv3)}\; f=x^{-6}y^{3}(1-x-y)  \\

 &  {\rm (lv4)}\; f=x^{-4}y^{2}(1-x-y)  & \\
& {\rm (lv5)}\; f=x^{-6}y^{3}(1-x-y)^{2} & \\

\end{array}$$
or 
$$\begin{array}{ll}
  (iii) f=x^{\frac{3}{k}}y^{\frac{l}{k}}(1+y), 
&(iv)  f=x^{-1+\frac{3}{k}-\frac{l}{k}}y^{\frac{l}{k}}(x+y) \\
\textrm{with} \ k \in \mathbb{Z}^{*}\setminus 3\mathbb{Z} \ \textrm{and}
\ l-k \in 3 \mathbb{Z},& \\

    (v) f=x^{\frac{4}{k}}y^{\frac{l}{k}}(1+y),  
&(vi) f=x^{-1+\frac{4}{k}-\frac{l}{k}}y^{\frac{l}{k}}(x+y) \\
\textrm{with}\  k \in \mathbb{Z}^{*}\setminus 2\mathbb{Z} \ \textrm{and}
\ l-k \in 4\mathbb{Z}& \\

    (vii) f=x^{\frac{6}{k}}y^{\frac{l}{k}}(1+y),
& (viii) f=x^{-1+\frac{6}{k}-\frac{l}{k}}y^{\frac{l}{k}}(x+y) \\
\textrm{with} \  k \in \mathbb{Z}^{*}, l \in 2\mathbb{Z} \ \textrm{and} \
kl-2 \in 6\mathbb{Z} & \\
(ix)  f=x^{\frac{4}{k}}y^{\frac{2l}{k}}(1+y)&
(x) f=x^{-1+\frac{4}{k}-\frac{2l}{k}}y^{\frac{2l}{k}}(1+y)\\
\textrm{with}\  k,l \in \mathbb{Z}^{*}\setminus 2\mathbb{Z}\\
(xi) f=x^{\frac{6}{k}}y^{\frac{3l}{k}}(1+y)&
(xii) f=x^{-1+\frac{6}{k}-\frac{3l}{k}}y^{\frac{3l}{k}}(1+y)\\
\textrm{with} \  k \in \mathbb{Z}^{*} \setminus 3\mathbb{Z} \ \textrm{and} \
l \in \mathbb{Z}^* \setminus 2\mathbb{Z} & \\

   \end{array}$$

and moreover $\rm gcd(k,l)=1.$

\end{montheo}
\textbf{Proof.}
An algebraic first integral is given by
$$f=x^{p_1}y^{p_2}(ax+by+c)^{r}, p_1, p_2 \in \mathbb{Z}, r \in \mathbb{N}^{*}.$$\\
This defines a divisor in $\mathbb{P}^2$: $D=p_1L_1+p_2L_2+rL_3$ where $L_i, i=1..
3$ are projective lines. As in previous section, the study below will depend on the geometry of the reduced divisor $\tilde{D}$ (i.e without multiplicities) associated to $D$. First we will consider the generic case, that is the projective lines $L_i, i=1..3$ have normal crossings toward each other.

\subsubsection{The case $a \neq 0, \ b \neq 0, \ c \neq 0.$}
First we may suppose under affine transformation $b=c=-a=1$.\\
The expression of $\tilde{D}$ invites us to divide the study in $4$:

\begin{equation}
\label{sit1}
p_1 >0, p_2 >0
\end{equation}

\begin{equation}
\label{sit2}
p_1 <0, p_2 >0, p_1+p_2+q >0
\end{equation}

\begin{equation}
\label{sit3}
p_1 <0, p_2 >0, p_1+p_2+q <0
\end{equation}

\begin{equation}
\label{sit4}
p_1 <0, p_2 >0, p_1+p_2+q =0.
\end{equation}

In the shape of (\ref{sit3}) the generic leaf is birational to the algebraic curve $X^{p_1}Y^{p_2}=t$ which is rational. Hence the generic case will be an obvious consequence of the $3$ following propositions:

\begin{prop}.
\label{prop1}
The algebraic curve  $x^{p}y^{q}(1-x-y)^r=1$
with $0 \leq p \leq q \leq r , \ \gcd(p,q,r)=1$ is of genus one 
if and only if $(p,q,r)=(1,1,1)$ or $(1,1,2)$ or $(1,2,3)$.
\end{prop}

\begin{prop}
\label{prop2}
The algebraic curve  $y^{q}(1-x-y)^r=x^{p}$
with $p,q,r >0$ , $-p+q+r < 0$ \  $\gcd(p,q,r)=1$ is of genus one if and only if $(p,q,r)=(1,2,2)$ or $(3,2,2)$.
\end{prop}

\begin{prop}
\label{prop3}
The algebraic curve  $y^{q}(1-x-y)^r=x^{p}$
with $p >0, q >r >0$ , $-p+q+r >0$ \  $\gcd(p,q,r)=1$ is elliptic if and only if $(p,q,r)=(3,1,1), (4,1,1), (4,2,1), (6,2,1), (6,3,1)$ or $(6,3,2)$.
\end{prop}




\textbf{Proof of Proposition \ref{prop1}.}\\
Let $\omega$ be a one form on a compact Riemann surface $S$. We write $\omega
= \sum_i a_iP_i$ with $P_i$ points of $S$. This sum is finite and we define the \textit{degree} of
$\omega$ : $\rm deg(\omega)=\sum_i a_i$.  According to the \textit{Poincar\'e-Hopf} formula (see \cite{algebraic}), any 1-form $\omega$ on
$S$ satisfies:
\begin{eqnarray}
\label{PH} \deg(\omega)=2g-2.
\end{eqnarray}
\\
\\
Now,  we use this formula with the riemann surface $\tilde{C}$ obtained after  desingularisation of the irreducible algebraic curve $C$ defined by the equation $x^{p}y^{q}(1-x-y)^{r}=1$. Let $ \pi : \tilde{C} \rightarrow C$ be such a desingularisation map. We compute below the degree of the one-form
 $\pi ^{*}\omega$  where (by abuse of notation):
$$\omega= -\frac{dx}{x[q-qx-(q+r)y]}=\frac{dy}{y[p-py-(p+r)x]}.$$ 

The $1$-form above has been chosen such that it has nor zeros nor poles outside the singular locus of $C$. Yet, $C$ is only singular in the three singular points meeting the line at infinity: $ [1:0:0], [0:1:0], [1:-1:0] .$
\\
First we investigate the local behavior of $\omega$ near $[0:1:0]$. We get local coordinates near  $ [1:0:0] $ as follows:
\\
Write $ x=\frac{1}{u} \ \ \textrm{with} \ \ u \rightarrow 0. $
\\
After this change of coordinates, the equation becomes:
$$y^{q}(u-1-yu)^{r}=u^{p+r}.$$
Since $u \rightarrow 0 $, we have the $m= gcd(q,p+r)$ different parametrisations of the $m$ local branches near this point:

$$ u=t^{\frac{q}{m}}$$
$$y=-e^{\frac{2ik \pi}{m}}t^{\frac{p+r}{m}}(1+o(t^{\frac{p+r}{m}})). \ \ k=0...m-1 .$$

For each branch, locally, $$\omega=-\frac{q}{m}t^{\frac{q}{m}-1}(1+o(t^{\frac{q}{m}-1}))dt.$$
Finally, for $\pi ^{*}\omega$ we get after a finite number of blowing-ups $m$ points where our $1$-form has a zero of order $\frac{q}{m}-1$.
\\
\\
The study is completely similar for the remaining singular points: near $[0:1:0]$ we obtain $n=gcd(p,q+r)$ points where the $1$-form $\pi ^{*}\omega $ has a zero of order $\frac{p}{n}-1$ and near $[1:-1:0]$, we have $l=gcd(r,p+q)$ points where $ \pi ^{*}\omega$ has a zero of order $\frac{r}{l}-1$.
\\
\\
Finally, the numbers involved satisfy the following relation:
\begin{eqnarray}
\label{dioph}
p+q+r-m-n-l=2g-2
\end{eqnarray}
and consequently, this curve is elliptic when:
\begin{eqnarray}
\label{dioph2}
p+q+r=m+n+l.
\end{eqnarray}

 Now we have to resolve this diophantine equation:
\\
\\
We always have:
$ m \leq q$, $n \leq p$ and $l \leq r .$
Hence (\ref{dioph2}) is true if and only if:
$$gcd(q,r+p)=q;$$
$$gcd(p,q+r)=p;$$
$$gcd(r,p+q)=r.$$
Let $\alpha , \ \beta , \ \gamma \in \mathbb{N}^{*}$ such that:
$$r+p=q \alpha; \eqno (a)$$
$$r+q=p \beta;  \eqno (b)$$
$$p+q=r \gamma.  \eqno (c)$$
Using $(a)$ and $(b)$, we obtain $(\alpha+1)q=(\beta+1)p$.
\\
Using $(b)$ and $(c)$, we obtain $(\gamma \beta-1)q=(\gamma+1)p$.
\\
\\
Hence we have:
$$\frac{q}{p}=\frac{\beta+1}{\alpha+1}=\frac{\gamma \beta-1}{\gamma+1}$$
which gives the following equation:
\begin{eqnarray}
\label{dioph3}
\alpha \beta   \gamma =2+\alpha + \beta +\gamma.
\end{eqnarray}
The solutions of this equation are under symmetry $(2,2,2), (3,3,1)$ and
$(5,2,1)$ which gives at last the solutions $(p,q,r)$ of Proposition \ref{prop1}.

\begin{flushright}
$\square$
\end{flushright}

\noindent \textbf{Proof of Proposition \ref{prop2}.}\\
The proof is similar. We still use
(\ref{PH}) with $\omega$ a $1-forme$ that without both zeros and poles outside the singular locus of the algebraic curve $C$ defined by $y^{q}(1-x-y)^{r}=x^{p}$:

$$\omega=
-\frac{dx}{x[q-qx-(q+r)y]}=\frac{dy}{y[-p+py-(r-p)x]}.$$ 

Here the singular locus is no longer as before. It has two singular points at infinity:
$[1:0:0]$ et $[1:-1:0]$ and moreover $(0,0)$ et $(0,1)$ in the affine chart. Local considerations as above naturally leads us to the following

\begin{lem}
The irreducible algebraic curve above has genus one if and only $p, \ q \ \textrm{and} \ r \ $ satisfy the equation:
\begin{eqnarray}
\label{dioph4}
q+r= gcd(p,q)+gcd(q,r)+gcd(r,r+q-p)+gcd(q,r+q-p).
\end{eqnarray}
\end{lem}

Writing $m=gcd(p,q)$, $n=gcd(q,r)$, $l=gcd(r,r+q-p)$ and $s=gcd(q,r+q-p)$, then there exists integers $\alpha, \beta, \gamma, \delta$ such that $r=n\alpha=l\beta$ and $q=s\gamma=m\delta$ so that $(\ref{dioph4})$  is equivalent to:
\begin{eqnarray}
\label{dioph5}
q+r=\frac{\gamma+\delta}{\gamma\delta}q+\frac{\alpha+\beta}{\alpha\beta}r
\end{eqnarray}
Hence we have $\gamma+\delta=\gamma\delta$ and $\alpha+\beta=\alpha\beta$ and therefore $\alpha=\beta=\gamma=\delta=2$.Hence $m=s$ and as $m$ divides $p$ and $q$ then $m$ divides $r$ and finally $m=1$. Similarly, we get $n=1$ and consequently $q=r=2$. Now, remind that $p < r+q=4$ so that $p$ either equals $1$ or $3$ ($2$ is excluded as $\gcd(p,q,r)=1$). Now we easily verify that $(1,2,2)$ and $(3,2,2)$ are the solutions to the equation $\ref{dioph4})$ above which proves Proposition \ref{prop2}.

\begin{flushright}
$\square$
\end{flushright}

\noindent \textbf{Proof of Proposition \ref{prop3}.}\\  
After the birational change of variable: $x \rightarrow \frac{1}{x}$ and $y \rightarrow \frac{y}{x}$, the genus (which is a birationnal invariant for curves) is the same as the genus of the algebraic curve:
$$x^{p-r-q}y^{q}(1-x-y)^{r}=1.$$ Then this is an immediate consequence of Proposition \ref{prop1} above.

\begin{flushright}
$\square$
\end{flushright}

Consequently we found $\rm (lv1-5)$ of Theorem \ref{th2}.

\subsubsection{The case $a=0.$}

\indent Under affine transformation we may suppose $b=c=1$. Geometrically
$\{y=0\}$ and $\{y=1\}$ both intersect  at infinity. Notice first the
following:  
$$x^{\lambda}y^{\mu}(1+y)=(xy^{n})^{\lambda}y^{\mu-n\lambda}(1+y)$$
 for $n \in \mathbb{Z}$, hence  $x^{\lambda}y^{\mu}(1+y)=t$ is birational to $x^{\lambda}y^{\mu-n\lambda}(1+y)=t$ so that we only need to study when $\lambda$ and $\mu$ are strictly positives. This naturally leads to the following: \\
\begin{prop}
\label{prop4}
The algebraic curve: $$C=\{(x,y) \in \mathbb{C}^{2},x^{p}y^{q}(1+y)^{r}=1\}$$ with $0\leq r \leq q $ and $ 0\leq p$,
where $gcd(p,q,r)=1.$ is elliptic if and only, under permutations of
$\{y=0\}$ and $\{y+1=0\}$ it is in the following list:\\
\\
$\begin{array}{ c c c c}
x^{3}y^{1+3u}(1+y)^{1+3v}=1; & x^{3}y^{2+3u}(1+y)^{2+3v}=1\\
x^{4}y^{1+4u}(1+y)^{1+4v}=1; & x^{4}y^{3+4u}(1+y)^{3+4v}=1; &
x^{4}y^{2(1+2u)}(1+y)^{r}=1, r \in Z^{*} \setminus 2 \mathbb{Z} \\
x^{6}y^{2+6u}(1+y)^{1+6v}=1; & x^{6}y^{5+6u}(1+y)^{4+6v}=1; &
x^{6}y^{3(2u+1)}(1+y)^{r}=1, r \in \mathbb{Z}^{*} \setminus 3\mathbb{Z}.
\end{array}$
\end{prop}
\textbf{Proof of Proposition \ref{prop4}.}\\
   We still use (\ref{PH}) with a judiciously chosen $\omega$ without zeros nor poles in its regular locus:
$$\omega=\dfrac{dx}{x(q+(q+r))y}=-\dfrac{dy}{py(1+y)}.$$
This curve has two points at infinity, namely $[1:0:0]$ and $[0:1:0]$, where $C$ is singular ($C$ is regular in the affine chart)
\\
Near $[1:0:0]$, we have two branches where a local equation of each is respectively:
$$Y^{q}=u^{p}$$
$$(1+Y)^{r}=u^{p}$$
where $x=\frac{1}{u}.$
Thus, writing: $m=pgcd(p,q)$ and $n=pgcd(p,r)$, we obtain the parametrisations:
$$Y=t^{\frac{p}{m}}$$
$$u=t^{\frac{q}{m}}$$
and
$$Y=t^{\frac{p}{n}}$$
$$u=t^{\frac{r}{n}}.$$
Both give a pole of order $1$ for $\pi^{*}\omega$, hence we obtain, adding up the different possible parametrisations, $-m-n$ in the Poincar\'e-Hopf formula.
\\
A similar calculus near the other point at infinity gives a zero of order $\frac{p}{l}-1$ where $l=(p,q+r)$
with $l$ different parametrisations.
\\
Thus we finally obtain the equality:
\begin{eqnarray}
\label{Dioph1}
p=m+n+l.
\end{eqnarray}

We want to resolve this equation.
Consider:
$$p=n\gamma$$
$$p=m\beta$$
$$p=l\alpha$$
with $\alpha, \beta, \gamma, \delta \in \mathbb{Z}^{*}$.
Then (\ref{Dioph1}) is equivalent to the following well-known equation:
\begin{eqnarray}
\label{Dioph2}
\alpha\beta\gamma=\alpha\beta+\alpha\gamma+\beta\gamma.
\end{eqnarray}
The solutions are up to permutation:
$$(3,3,3)$$
$$(2,4,4)$$
$$(2,3,6).$$
The solution $(3,3,3)$ implies $m=n=l.$
\\
As $\gcd(p,q,r)=1$ we have: $m=n=l=1$ and so $p=3$, $1=(r,3)$, $1=(q,3)$. Hence $r=1\ , \ 2 \ \bmod(3)$ and so does $q$. Finally, reminding $l=\gcd(p,q+r)$, we conclude  that  $(p,q,r)=(3,1+3u,1+3v),(3,2+3u,2+3v)$.
\\
\\
The solution $(2,4,4)$  implies  $l=2n$ and $m=n=1$ and the same
argument shows that $(p,q,r)=(4,1+4u,1+4v)$ or $(4,3+4u,3+4v).$ There
are $2$ other solutions (permutations of $(2,4,4)$). Under
permutations of the two lines $\{y=0\}$ and $\{y+1=0\}$, we only need
to study $(4,2,4)$. A similar resolution thus gives $(p,q,r)=(4,2(1+2u),1+2v).$
\\
\\
The solution $(2,3,6)$ implies $m=2$, $n=1$, $l=3$, so
$(p,q,r)=(6,2+6u,1+6v)$ or $(6,5+6u,4+6v).$ As above, we need to take
under consideration the solutions $(3,2,6)$ and $(6,2,3)$ wich
respectively gives $(p,q,r)=(6,3(2u+1),r)$ with $ { \rm gcd}(r,3)=1$
and $(6,2(3u+1),6v+1)$ or $(6,2(3u+2),6v+5)$. Finally the proposition is proved.
\begin{flushright}
$\square$
\end{flushright}
This we obtain the last cases of the left column of Theorem \ref{th2}.

\subsubsection{The case $c=0.$}
Under affine transformation we may suppose $a=b=1$.\\
Here the three lines $\{x=0\}$, $\{y=0\}$ and $\{x+y=0\}$ intersect themselves at the origine.\\
Now, the algebraic curve $x^{\lambda}y^{\mu}(x+y)=t $ is obviously birational to $x^{\lambda+\mu+1}y^{\mu}(1+y)=t$ , so this case falls from the preceding results and we get the last cases of the right column of Theorem \ref{th2}. \\
\\
We have investigated in fact all the possible first integrals. Indeed, if our foliation admits a first integral: $f=x^{-\alpha}y^{-\beta}(ax+by+c)$ with $\alpha, \beta$ real positive numbers, then after affine transformation it has a first integral: $g=X^{-\frac{1}{\alpha}}Y^{\frac{\beta}{\alpha}}(AX+BY+C)$. Hence Theorem \ref{th2} is proved.
\begin{flushright}
$\square$
\end{flushright}

\subsection{The reversible Lotka-Voltera case}
\label{section23}
\begin{montheo}
\label{theo}
\label{th3}
The reversible Lotka-Voltera foliation $\mathcal{F}(df)$ is elliptic if and only if, after an affine change of
the variables, it has a first integral of the form
$$
\begin{array}{ll}
\rm (rlv1)\; f=xy(1-x-y) & \rm(rlv2) \;f=x^{-3}y(1-x-y)\\
\rm(rlv3)\;f=x^{2}y(1-x-y)&
\rm(rlv4)\;f=x^{-4}y(1-x-y)\\
\rm(rlv5)\;f=x^{-3}y^{2}(1-x-y)^{2}&
\rm(rlv6)\;f=x^{-1}y^{2}(1-x-y)^{2}\\

\\
or:
\\
(xiii) f=x^{\frac{3}{k}}(y+1)(y-1), k\in \mathbb{Z}^*\setminus 3\mathbb{Z} &
(x1v) f=x^{-2+\frac{3}{k}}(x+y)(x-y), k\in \mathbb{Z}^*\setminus 3\mathbb{Z}\\
(xv) f=x^{\frac{4}{k}}(y+1)(y-1),
k\in \mathbb{Z}^*\setminus 2\mathbb{Z} &
(xvi) f=x^{-2+ \frac{4}{k}}(x+y)(x-y),
k\in \mathbb{Z}^*\setminus 2\mathbb{Z}\\
\end{array}
$$
\end{montheo}
\textbf{Proof.}\\
 Without loss of generality we may suppose that the foliation has a  first integral of the form
$f =x^{\lambda}(y^{2}+a(x-\frac{b}{2a})^{2})$ if $a \neq 0$ and $f=x^{\lambda}(y^{2}+c)$ otherwise (we deal about quadratic foliations so that $c$ necessary does not vanish). For $a =0$ this a consequence of Proposition \ref{prop4}. Now look at $a \neq 0$ :\\
\\
if $b=0$, the  curve $x^{\lambda}(y^{2}+ax^{2})=t$ is  birational to $x^{\lambda+2}(y^{2}+1)=t$ thus the conditions are  $p+2q=3,4$ or $-2q-p=3,4$.\\
if $b \neq 0$ the quadric is a reducible polynomial so that this case is a straightforward consequence of 
Propositions \ref{prop1}, \ref{prop2} and \ref{prop3}. We notice that the last cases of reversible Lotka-Volterra  are exactly the last cases of Lotka-Volterra under the condition $l=k$ when it is possible (For $(vii)$ and $(viii)$ of Theorem \ref{th2} we can't have $k=l$). This gives Theorem \ref{th3}.
\begin{flushright}
$\square$
\end{flushright}
Notice that the case $b=0$ is also a consequence of the  degenerate Lotka-Voltera case with the three invariant lines involved intersecting themselves, but the calculus is here so easy that we proved it directly and is useful to test our preceeding calculus.\\

\section{Topology of the singular fibers and Kodaira's classification}
\label{sec3}
Now we focus on the singular fibers of the induced elliptic surfaces. First
of all, Recall that two birational elliptic surfaces have the same
minimal model (see \cite{K75,mir}). Some of our previous elliptic surfaces are obviously birationnals and therefore have the same singular fibers under permutation. First we investigate such mappings. Then we give some examples of computation of the singular fibers to illustrate the way we obtained Tables \ref{tab2}, \ref{tab3}, \ref{tab4}, \ref{tab5}, \ref{tab6}, \ref{tab7}, \ref{tab8}.
\subsection{The reversible case}
\subsubsection{Birational mappings}
The first integrals are given by the algebraic equation:
$$x^{\lambda}(y^{2}+ax^{2}+bx+c)=t$$ with $a$, $b$, $c$ complex numbers satisfying some conditions and $\lambda$ a rational number.
We have an easy birational mapping (we already used it, see Section \ref{213}): $X=\frac{1}{x}$, $Y=\frac{y}{x}$ which leads to $$x^{-2-\lambda}(y^{2}+cx^{2}+bx+a)=t.$$
When considering this mapping in $\mathbb{P}^2$ with homogeneous
coordinates $[x:y:z]$ this last permutes in fact the projactive lines
$\{x=0\}$ and $\{z=0$\}. Thus for each line of Table \ref{tab1} we
only need to study either the right or the left element.\\
\\
For degenerate cases, notice the change of variables $(X,Y)=(xy,y)$
birationnally leads $(i)$ (resp. $(ii)$) to $(x)$ (resp. $(xii)$)
Lotka-Volterra elliptic case with $l=1$. Consequently, such cases will
be a consequence of the calculus of the singular fibers of the
Lotka-Volterra cases (see below). 
\begin{table}
\begin{center}
\begin{tabular}{l | l}
$\lambda <-1$ &$\lambda>-1$\\
\hline
\\
$(rv1) \ x^{-3}(y^{2}+ax^{2}+bx+c)=t$&
$(rv2)  \ x(y^{2}+cx^{2}+bx+a)=t$\\
\\
$(rv3) \ x^{-\frac{3}{2}}(y^{2}+ax^{2}+bx+c)=t, c \neq 0$&
$(rv4) \ x^{-\frac{1}{2}}(y^{2}+cx^{2}+bx+a)=t, c \neq 0$\\
\\
$(rv5)  \ x^{-4}(y\index{}^{2}+ax^{2}+bx+c)=t$&
$(rv6)  \ x^{2}(y^{2}+cx^{2}+bx+a)=t$\\
\\
$(rv7)  \ x^{-\frac{4}{3}}(y^{2}+bx+c)=t$&
$(rv8) \ x^{-\frac{2}{3}}(y^{2}+cx^{2}+bx)=t$\\
\\
$(rv9)  \ x^{-\frac{4}{3}}(y^{2}+ax^{2}+bx)=t$&
$(rv10)  \ x^{-\frac{2}{3}}(y^{2}+bx+a)=t$\\
\\
$(rv11) \ x^{-\frac{5}{3}}(y^{2}+ax^{2}+bx)=t$&
$(rv12) \ x^{-\frac{1}{3}}(y^{2}+bx+a)=t$\\
\\
$(rv13)  \ x^{-\frac{5}{4}}(y^{2}+ax^{2}+bx)=t$&
$(rv14)  \ x^{-\frac{3}{4}}(y^{2}+bx+a)=t$\\
\\
$(rv15) \ x^{-\frac{7}{4}}(y^{2}+ax^{2}+bx)=t$&
$(rv16) \ x^{-\frac{1}{4}}(y^{2}+bx+a)=t$\\
\\
$(rv17) \ x^{-\frac{5}{2}}(y^{2}+ax^{2}+bx)=t$&
$(rv18) \ x^{\frac{1}{2}}(y^{2}+bx+a)=t$\\
\\
\hline
\\

\end{tabular}
\begin{tabular}{l l}
$ (i) \ x^{-1+ \frac{2}{k}}(y^{2}+x)=t$ &
$k\in\mathbb{Z}^*\setminus 2\mathbb{Z}$\\
\\
$(ii) \ x^{-1+ \frac{3}{k}}(y^{2}+x)=t$&
$k\in\mathbb{Z}^*\setminus 3\mathbb{Z}$\\
\\
\hline
\end{tabular}
\end{center}
\caption{\label{tab1}The elliptic reversible cases.}
\end{table}
\vspace{0.5cm}

Be careful that the geometry of the divisors appearing in the first integrals (including the line at infinity) is of importance as we shall blow-up the indetermination points. Birationally, the  different geometrical description of the divisors in the reversible case are the following:\\
\\
$(1)$ The divisors are in general position (see Figure \ref{Config}). This concerns $\rm(rv2), \rm(rv4), \rm(rv6)$ with $a, b, c \neq 0$.\\
\\
$(2)$  $\{Q=0\}$ and $\{x=0\}$ are in general position and $\{Q=0\}$ and $\{z=0\}$ have only one tangent double point (see Figure \ref{example2}). This concerns $\rm(rv2)$ with $c=0$, $\rm(rv3)$ with $a=0$, $\rm(rv6)$ with $c=0$, $\rm(rv7), \rm(rv10), \rm(rv12), \rm(rv14), \rm(rv16)$.\\
\\
$(3)$ Both projective lines $\{x=0\}$ and $\{z=0\}$ have a double tangent point with the quadric . This concerns $\rm(i)$ and $\rm(ii)$. 
\subsubsection{The singular fibers $t=0$ and $t=\infty$}

\begin{figure}[h]

\begin{center}

\input{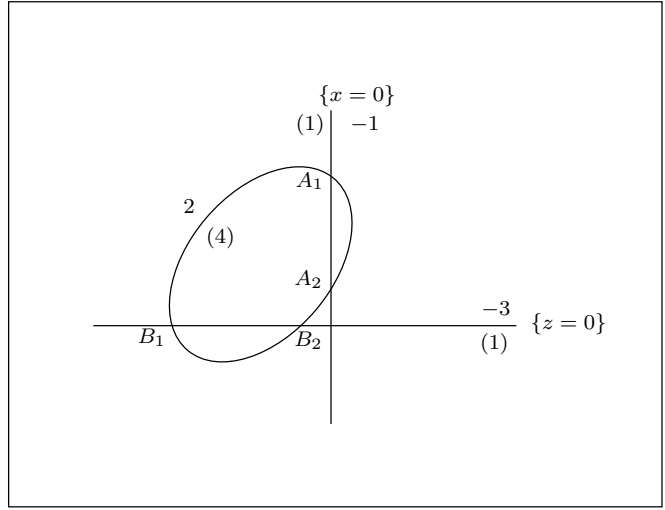} 
\caption{\label{Config}Geometrical situation of $\rm(rv4)$.}
\end{center}

\end{figure}

Here we illustrate the results with examples:
\\
\\
\textit{EXAMPLE 1: The singular fibers of $\rm(rv4)$:}\\
\\
Embedding our first integral in $\mathbb{P}^{2}$, one have:
$$\dfrac{(y^{2}+ax^{2}+bxz+cz^{2})^{2}}{xz^{3}}=t.$$
Geometrically, there are two lines: $\{x=0\}$ and $\{z=0\}$ with  multiplicities respectively $-1$ and $-3$ that intersect in $[0:1:0]$, and a conic that intersects both lines in four points, namely $A_1=[ 0:\sqrt{-c}:1 ] $, $A_2=[ 0:-\sqrt{-c}:1 ]$, $B_1=[\sqrt{-a}:1:0]$, $B_2=[-\sqrt{-a}:1:0]$ with normal crossing each time (see Figure \ref{Config}).\\

 The rational function is not defined in these four points, thus we need to blow-up them.
\begin{itemize}
 \item Near $A_1$
\\
In local coordinates the rational function becomes: $\frac{Y^{2}}{X}$.\\
We need two blowing-ups to define the rational function near this point:
$$\frac{Y^{2}}{X}\rightarrow \frac{Y}{X} \rightarrow \textrm{separation of both local branches}$$
Remind that blowing a point of  $\mathbb{P}^{2}$ that belongs to a divisor $D$ decreases the self-intersection of $D$ by one (see \cite{algebraic} for example). Writing the self-intersection and the multiplicities (the self-intersection numbers are inside the () ) we get the situation of Figure \ref{fig3}.\\

\begin{figure}[h!]
\input{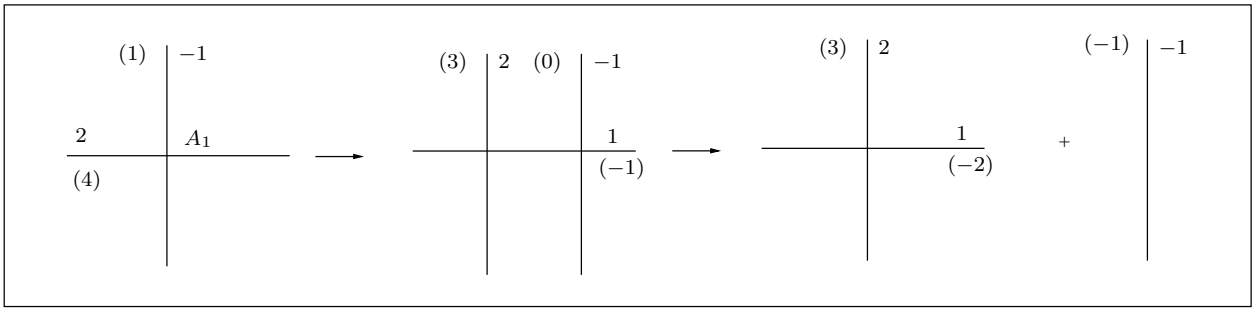}
\caption{\label{fig3}The successive blowing-ups of $A_1$.}
\end{figure}
We study $A_2$ along the same lines. See Figure \ref{fig4}.\\

\begin{figure}[h!]
\begin{center}
\input{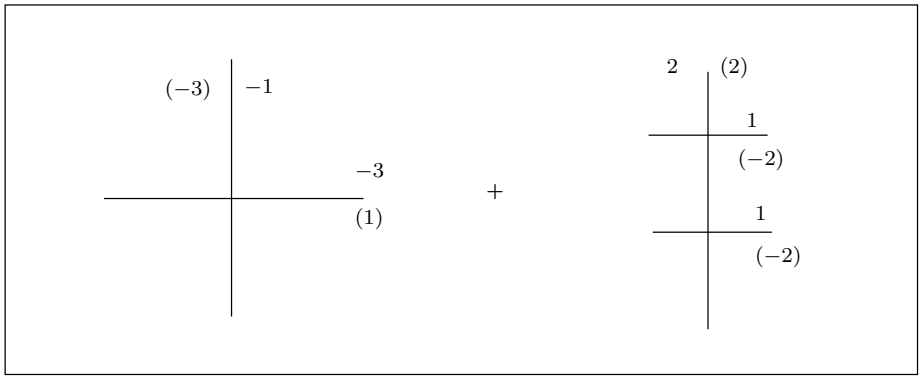}
\caption{\label{fig4}Summary ot the situation after blowing-up $A_1$ and $A_2$.}
\end{center}
\end{figure}
\item Near $B_1$.
\\
 Locally the rational function becomes $\frac{Y^{2}}{Z^{3}}.$\\
The successive blowing-ups give the following local equations until separation:
$$\frac{Y^{2}}{Z^{3}} \rightarrow  \frac{Y^{2}}{Z} \rightarrow \frac{Y}{Z} \rightarrow \textrm{separation of the branches }$$
see Figure \ref{fig5}. \\

\begin{figure}[h!]
\begin{center}
\input{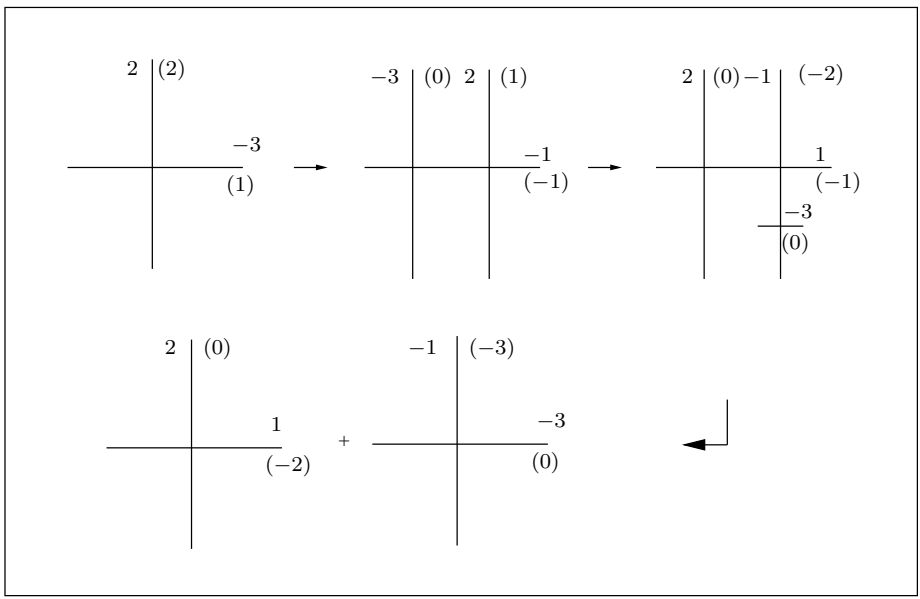}
\caption{\label{fig5}The successive blowing-ups of $B_1$.}
\end{center}
\end{figure}
The situation is still the same near $B_2$. We obtain $2$ singular fibers. See Figure \ref{fig6}. \\
\\


\begin{figure}[h]
\begin{center}

\input{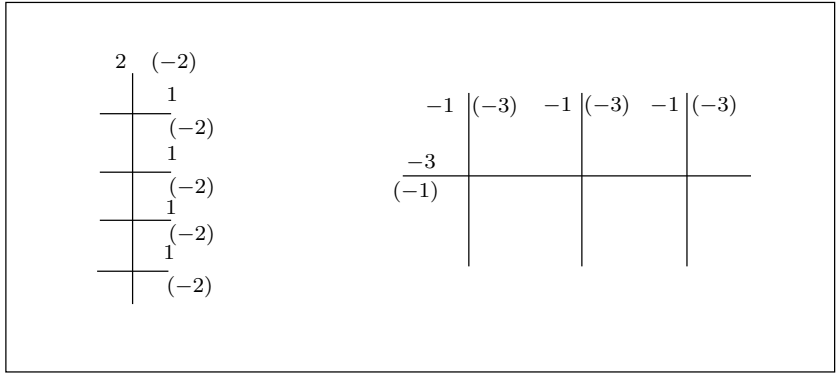}
\caption{\label{fig6}The fibers $t=0$ and $t= \infty$ of $\rm(rv4)$.}
\end{center}
\end{figure}

\end{itemize}

\begin{figure}[h!]
\begin{center}
\input{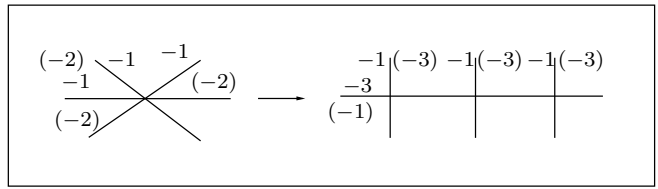}
\caption{\label{fig7}Contraction of the fiber $t=\infty$ of $\rm(rv4)$.}
\end{center}
\end{figure}

We now have to recognize these singular fibers in Kodaira's classification (see \cite{singularbis}). The singular fiber  $t=0$ is $I_0^{*}$, but we don't recognize the other. This is because we have branches with self-intersection $-1$. Recall that Kodaira's classification involves \textit{minimal} elliptic surfaces i.e no fiber contains an exceptionnal curve of the first kind. \\Finally we get $IV$ (the singular fiber at infinity). See Figure \ref{fig7}.\\

\begin{figure}[h!]
\begin{center}
\input{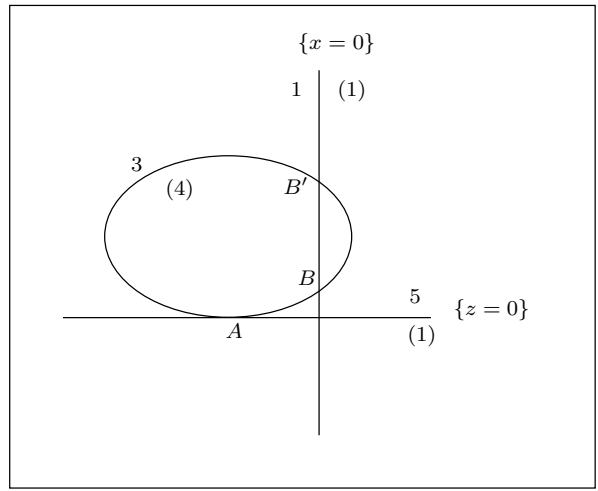}
\end{center}
\caption{\label{example2}The divisor associated to $\rm (rv12)$.}

\end{figure}
\textit{EXAMPLE 2 : The singular fibers of $\rm(rv12)$ :}\\
\\
The rational funcion here is:
$$\frac{(y^{2}+bxz+az^{2})^{3}}{xz^{5}}.$$
The geometrical situation is explained in Figure \ref{example2}.

\begin{figure}

\begin{center}

\includegraphics[height=12cm, width=14cm]{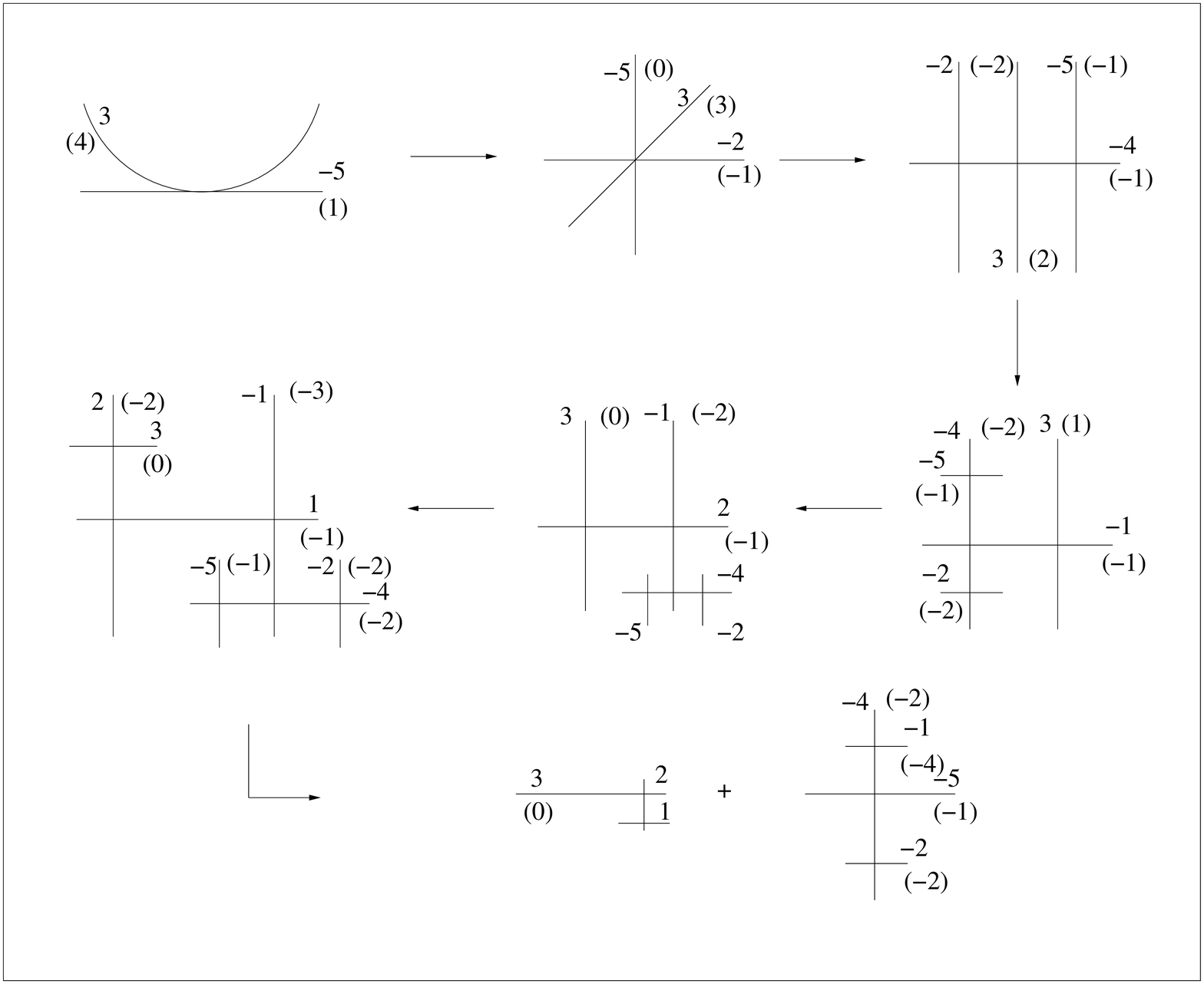}\\

\caption{\label{blowex2}Successive blowing-ups of $\rm(rv12)$ near $A$.}
\end{center}

\end{figure}

\begin{figure}
\begin{center}
\includegraphics[height=4cm, width=6cm]{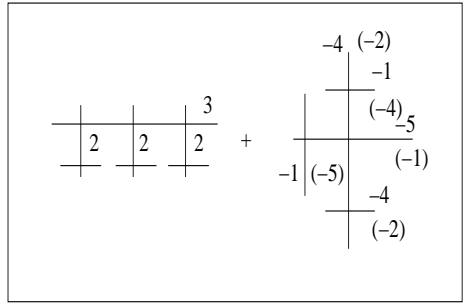}
\end{center}
\caption{\label{fib2}The fibers at $t=0$ and $t= \infty$ of $\rm(rv12)$.}

\end{figure}

\begin{figure}

\includegraphics[height=4cm, width=15cm]{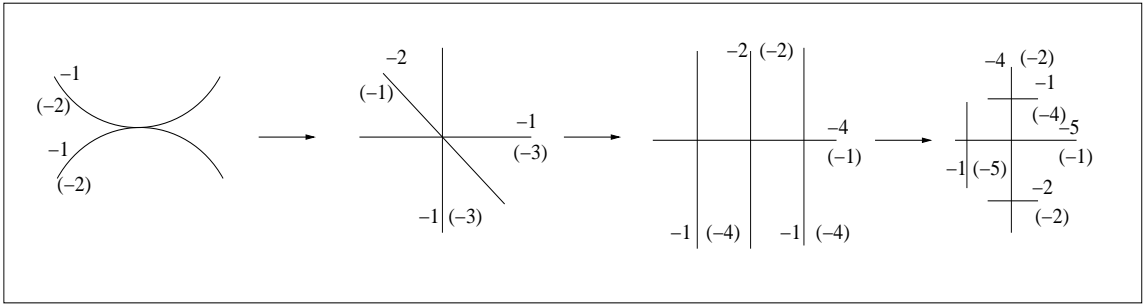}

\caption{\label{con2}Contraction of divisors with self-intersection $-1$ for the fiber at infinity of $\rm(rv12)$.}

\end{figure}

\begin{itemize}
   \item Study near $A$:\\
Locally the rational function becomes: $\frac{Y^{3}}{Z^{5}}.$
To begin with, we have:
$$ \frac{(Y^{2}+Z)^{3}}{Z^{5}} \rightarrow \frac{(Y+Z)^{3}}{Z^{5}Y^{2}} \rightarrow  \frac{(Z+1)^{3}}{Z^{5}Y^{4}} \rightarrow \textrm{separation of both local branches} $$

Next we have to blow-up  the point with local coordinate: $(Y=0,Z=-1)$, what gives  locally :
$$\frac{Z^{3}}{Y^{4}} \rightarrow \frac{Z^{3}}{Y} \rightarrow \frac{Z^{2}}{Y} \rightarrow \frac{Z}{Y} \rightarrow \textrm{separation of both local branches}.$$
the geometrical explanations  are given in Figure \ref{blowex2}.

   \item Study near $B$:\\
Locally the rational function becomes: $\frac{Y^{3}}{X}$ . Such a calculus has already been done in Example 1.

\end{itemize}
We finally obtain two fibers (see Figure \ref{fib2}).
\\
For $t=0$ we recognize $IV*$. For $t=\infty$, one have to contract divisors with self-intersection $-1$ as in Figure \ref{con2}. We finally get $III$ of Kodaira's classification.
\subsubsection{Other(s) singular fiber(s)}
Be careful that we do not have the full list of singular fibers: we also have to consider the singular points of our foliation whiches do not intersect both lines and the conic curve above. Here, the results concern the whole class of reversible systems inducing elliptic fibrations:\\
Writing: $$f=x^{\lambda}(y^{2}+ax^{2}+bx+c),$$
$(x,y)$ is a singular point if and only : $$x^{\lambda}y=0$$
$$x^{\lambda-1}(\lambda y^{2}+(\lambda +2)ax^{2}+(\lambda+1)bx+\lambda c)=0.$$
Consequently $(x_0,y_0)$ is a singular point which does not intersect both lines and the conic curve above if and only:
$$y_0=0$$
$$ x_0 \  \textrm{is a zero of the polynomial }: \ P=(\lambda +2)ax^{2}+(\lambda+1)bx+\lambda c \ \textrm{and} \ x_0 \neq 0.$$
Let $\delta$ be the discriminant of $P$. We resume the general happening  below:\\
\\
If $a \neq 0$, reminding $P(x_0)=0$, we get:\\
\\
$\frac{\partial^{2}{f}}{\partial{x^{2}}}(x_0,0)=x_0^{\lambda-1}P'(x_0)$\\
$\frac{\partial^{3}{f}}{\partial{x^{3}}}(x_0,0)=x_0^{\lambda-2}(2(\lambda-1)P'(x_0)+2a(\lambda+2)x_0 ).$\\
\\
As $\lambda \neq 0,\ -1 $ in the lists obtained concerning the whole reversible case above, different critical points have different critical values and  $P'(x_0)=0 \Leftrightarrow \delta =0 \Leftrightarrow (\lambda+1)b^{2}=4\lambda ac$, we have  the following (remind that for the reversible case $b^{2}-4ac \neq 0$):
\begin{lem} For $a\neq0$ and $c \neq 0$, if $\delta \neq 0$  we obtain two different singular curves with
a normal crossing, that is $I_1$ in Kodaira's classification and if $\delta=0$, we obtain one singular fiber with a cusp, that is $II$. If $a \neq 0$ and $c=0$, or $a=0$ and $c \neq 0$, we obtain one singular fiber with normal crossing, that is $I_1$. Otherwise, there are no more singular fibers.
\end{lem}

Finally, we are now able to compute all the singular fibers. The results are given in Tables \ref{tab2}, \ref{tab3} and \ref{tab4}.


\begin{table}[h!]
\begin{center}
\begin{tabular}{| l | c | c | c | c |}
\hline
Fibration & $\{t=\infty\}$ & $\{t=0\}$ & $t_1$ & $t_2$\\
\hline
$\rm(rv1)$ & $IV*$ & $I_2$ & $I_1 $ & $I_1$ \\
$\rm(rv2)$ & $IV*$ & $I_2$ & $I_1 $ & $I_1$\\
\hline
$\rm(rv3)$ & $IV$ & $I_0*$ & $I_1 $ & $I_1$\\
$\rm(rv4)$ & $IV$ & $I_0*$ & $I_1 $ & $I_1$\\
\hline
$\rm(rv5)$ & $III*$ & $I_1$ & $I_1 $ & $I_1$ \\
$\rm(rv6)$ & $III*$ & $I_1$ & $I_1 $ & $I_1$\\
\hline
\end{tabular}
\caption{\label{tab2}The elliptic reversible case with $4$ singular fibers. }
\end{center}
\end{table}

\begin{table}[h!]
\begin{center}
\begin{tabular}{| c | c | c | c | }
\hline
Fibration & $\{t=\infty\}$ &  $\{t=0\}$ & $t_1$ \\
\hline
$\rm(rv1)$  $a, b, c \neq  0$,  $\delta =0$ & $IV*$  & $I_2$ &  $II$ \\
 \ \ \ \ \ \ \                 $a \neq 0 $ ,$c=0$ & $III*$ & $I_2$ & $I_1$ \\
     \ \ \ \ \ \ \                          $a=0$, $c \neq 0$ &  $IV*$   &     $III$     &  $I_1$         \\
$\rm(rv2)$ $a, b, c \neq  0$, $\delta =0$ & $IV*$  & $I_2$ &  $II$ \\
      \ \ \ \ \ \ \                         $a \neq 0 $ ,$c=0$ & $III*$ & $I_2$ & $ I_1$ \\
 \ \ \ \ \ \ \            $a=0$, $c \neq 0$ & $IV*$       &   $III$       &   $I_1$     \\
\hline
$\rm(rv3)$  $a, b, c \neq  0$, $\delta =0$ & $IV$ & $I_0*$    & $II$ \\
          \ \ \ \ \ \ \                    $a=0$                                  &   $IV$         &  $I_1*$       &    $I_1$    \\
$\rm(rv4)$, $a, b, c \neq  0$, $\delta =0$ & $IV$ & $I_0*$ & $II$ \\
            \ \ \ \ \ \ \                   $a=0$ & $IV$ & $I_1*$ &    $I_1$                                                            \\
\hline
$\rm(rv5)$  $a, b, c \neq  0$, $\delta =0$ & $III*$ & $I_1$ &  $II$ \\
             \ \ \ \ \ \ \                  $a \neq0$ ,$c=0$ & $II*$ &$I_1$  & $I_1$ \\
               \ \ \ \ \ \ \                $a=0$, $ c \neq 0$ &  $III*$   &     $II$     &  $I_1$         \\

$\rm(rv6)$  $a, b, c \neq  0$, $\delta =0$ & $III*$ & $I_1$ &  $II$ \\
  \ \ \ \ \ \ \            $a \neq0$ ,$c=0$ & $II*$ &$I_1$  & $I_1$ \\
     \ \ \ \ \ \ \                          $a=0$, $ c \neq 0$ & $III*$    &      $II$    &    $I_1$       \\
\hline
$\rm(rv7)$ & $IV*$ & $III$ &  $ I_1 $ \\
 $\rm(rv8)$ & $IV*$ & $III$ &  $ I_1 $ \\
\hline
$\rm(rv9)$ & $I_1*$ & $IV$ &  $ I_1  $ \\
 $\rm(rv10)$ & $I_1*$ & $IV$ &  $ I_1  $ \\
\hline
$\rm(rv11)$ & $III$ & $IV*$ &  $ I_1  $\\
 $\rm(rv12)$ & $III$ & $IV*$ &  $ I_1  $\\
\hline
$\rm(rv13)$ & $IV*$ & $III$ &  $ I_1 $\\
 $\rm(rv14)$ & $IV*$ & $III$ &  $ I_1 $\\
\hline
$\rm(rv15)$ & $III*$ & $II$ &  $ I_1  $\\
 $\rm(rv16)$ & $III*$ & $II$ &  $ I_1  $\\
\hline
$\rm(rv17)$ & $II*$ & $I_1$ &  $ I_1  $.\\
$\rm(rv18)$ & $II*$ & $I_1$ &  $ I_1  $.\\
\hline
\end{tabular}
\end{center}
\caption{\label{tab3}The elliptic reversible case with $3$ singular fibers.}
\end{table}
\begin{table}[h!]
\begin{center}
\begin{tabular}{| c | c | c |}
\hline
Fibration & $\{t=\infty\}$ &  $\{t=0\}$  \\
\hline
$(i), k=1 \mod(4)$ & $III^*$  & $III$ \\
\ \ \ \ \ $k=3 \mod(4)$ & $III$ & $III^*$\\
\hline
$(ii)$ $k=1,2 \mod(6)$ &  $II^*$ & $II$   \\
\ \ \ \ \  $k=3,5 \mod(6)$  & $II$ & $II^*$  \\
\hline
\end{tabular}
\end{center}

\caption{\label{tab4}The reversible case with $2$ singular fibers.}
\end{table}
\newpage
\subsection{The Lotka-Volterra case}
\subsubsection{Birational mappings}
The same birational mapping: $X=\frac{1}{x}$ and $Y=\frac{y}{x}$,
leads:
$$x^{\lambda}y^{\mu}(ax+by+c)=t$$ birationally to:
$$x^{-\lambda-\mu-1}y^{\mu}(cx+by+a)=t.$$
Using this, one can immediately verify that $\rm(lv1)$, $\rm(lv2)$, $\rm(lv3)$ and $\rm(lv5)$ are birationals and $\rm(lv4)$ is birational to $\rm(rlv3)$ and $\rm(rlv4)$.
For the last Lotka-Volterra cases, there are another obvious birational mappings:
$$X=xy^{ u}(1+y)^{v} \ , \ Y=y$$
or:
$$Y=xy^{ u}(1+y)^{ v} \ , \ X=x.$$
with $u, v \in \mathbb{Z}$ judiciously chosen.
\\

\subsubsection{The fibers $t=0$ and $t={\infty}$}
\begin{figure}[h]
\begin{center}
\input{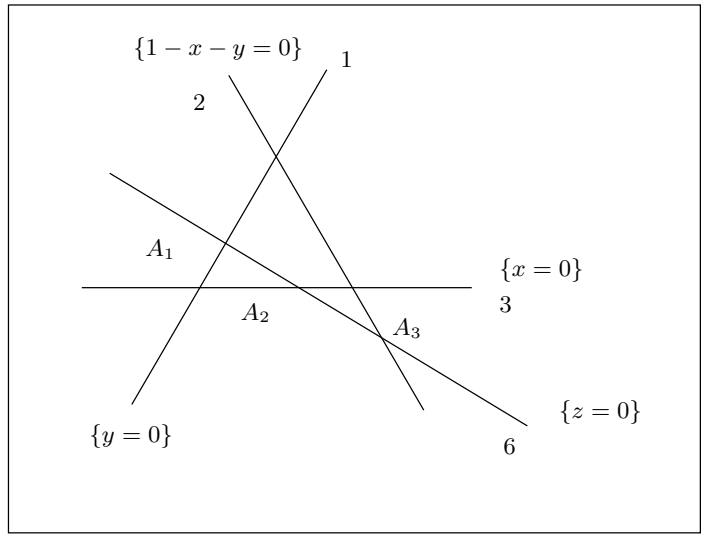}
\caption{\label{glv}The geometrical situation of $\rm(lv1)$.}
\end{center}
\end{figure}

\begin{figure}[h]
\begin{center}
\includegraphics[height=3cm,width=7cm]{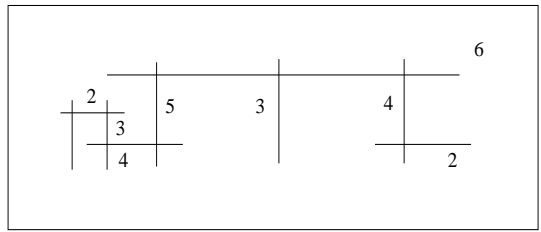}
\caption{\label{IIstar}The fiber at infinity for $\rm(lv1)$.}
\end{center}
\end{figure}

\begin{figure}[h!]
\begin{center}
\includegraphics[height=4.5cm,width=5.5cm]{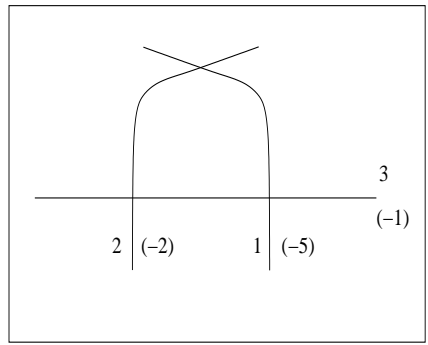}
\caption{\label{fib}The fiber  $t=0$ for $\rm(lv1)$.}
\end{center}
\end{figure}

\begin{figure}[h!]
\begin{center}
\includegraphics[height=10cm,width=5.5cm]{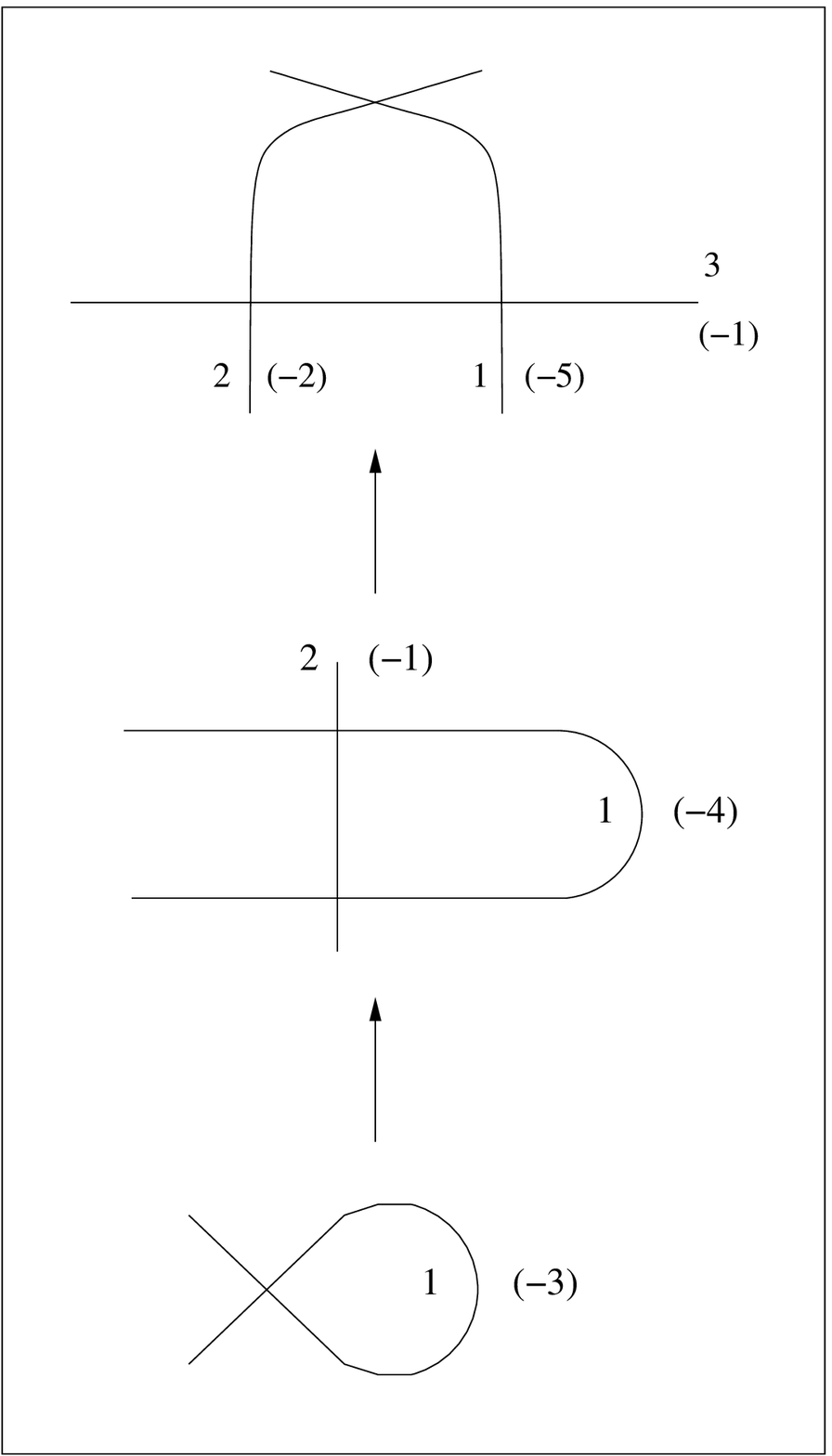}
\caption{\label{contraction}Contraction of the fiber $t=0$ for $\rm(lv1)$.}
\end{center}
\end{figure}

\textit{EXAMPLE 3:The singular fibers of $\rm(lv1)$}
\\
Here the rational function is:
$$\frac{x^{2}y^{3}(z-x-y)}{z^{6}}.$$
The geometrical situation is explained in Figure \ref{glv}. The intersection points with opposite multiplicity need to be blown-up. Here there are $3$ points: $A_1=[0:1:0]$, $A_2=[1:0:0]$, $A_3=[1:-1:0]$.
\\
Near $A_1$, locally the rational function becomes $\frac{X^{2}}{Z^{6}}$ such that we only need three
blowing-ups to separate local branches. This situation is well-known like near $A_2$ and $A_3$, where
we need respectively $2$ and $6$ blowing-ups. See Figures \ref{IIstar}
and \ref{fib}. The fiber at infinity is $II^{*}$. For the fiber $t=0$, we need two contractions as explained in
Figure \ref{contraction} and we finally get $I_1$.
\\

\textit{EXAMPLE 4: The singular fibers of $(vii)$ with $l=4 \bmod (6)$.}
Under birational equivalence, the rational function we need to consider is:
$$\frac{y^{4}(z+y)^{5}}{x^{6}z^{3}}.$$

\begin{figure}[h!]
\begin{center}
\input{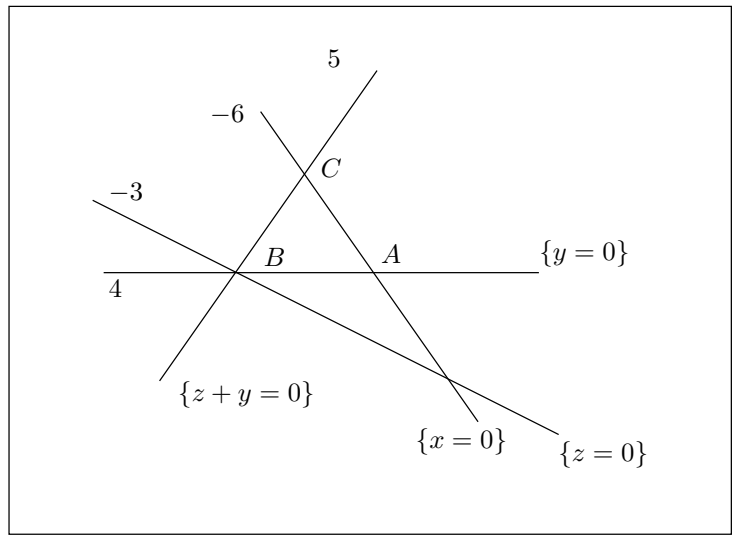}
\end{center}
\caption{\label{youpi}Geometrical situation of $(vii)$.}
\end{figure}
\begin{figure}[h!]
\begin{center}
\includegraphics[height=9cm,width=4.5cm]{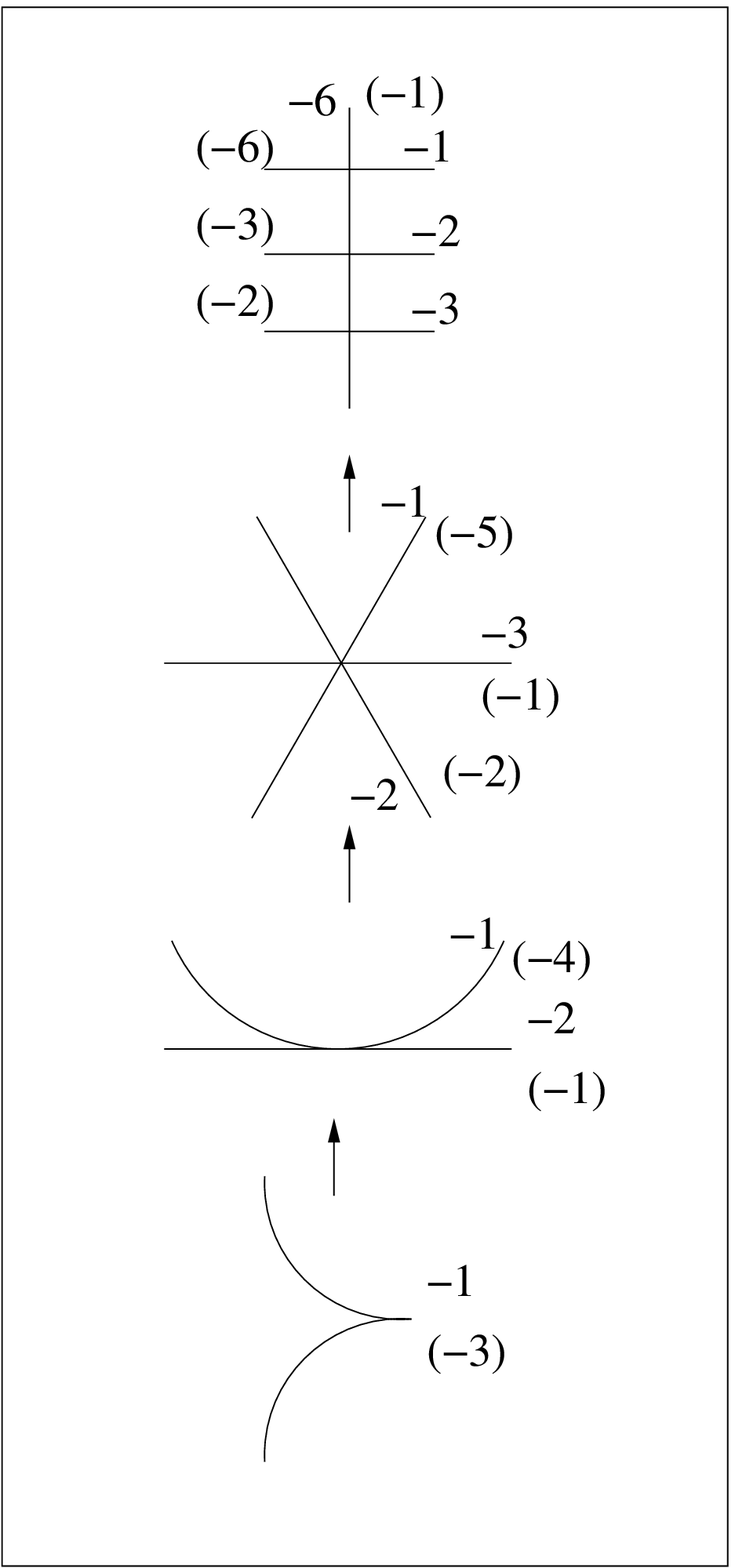}
\end{center}
\caption{\label{final}Contractions of $t=0$ for $(vii)$.}
\end{figure}

We have to blow up $3$ points: $A=[0:0:1]$, $B=[1:0:0]$ and $C=[0:1:-1].$
For $A$ and $C$ we have normal crossings and the situation is similar to precedent ones. We need to pay little more attention for the blowing-up of $B$:\\
Locally the rational function becomes:  $\frac{Y^{4}(Y+Z)^{5}}{z^{3}}$. Here the first blowing- up separates the three branches. Now we need to blow-up the intersection point of the branch with multiplicity $6$ and the branch with multiplicity $-3$. Locally the rational function is: $\frac{Y^{6}}{Z^{3}}$ and we get in a well-known situation. We obtain $II*$ for $t=0$. For $t=\infty$ we need three contractions to finally obtain $II$ (see Figure \ref{final}).
\subsubsection{Other singular fiber}
We write:
$$f=x^{\lambda}y^{\mu}(ax+by+c).$$
This is here elementary linear algebra and it immediately gives the following:
\begin{lem}
Under the assumptions $\lambda, \mu \neq 0, \lambda \neq -1, \mu \neq -1$ and $\lambda + \mu +1 \neq 0$,
the function above gives rise to another singular fiber if and only if $a, b, c \neq 0$ and the corresponding singular fiber is $I_1$.
\end{lem}
\begin{rem}
Such assumptions hold for our Lotka-Volterra and reversible Lotka-Volterra systems inducing elliptic fibrations.
\end{rem}
\begin{table}[h!]
\begin{center}
\begin{tabular}{| c | c | c | c |}
\hline
Fibration & Fiber $\{t=\infty\}$ & Fiber $t=0$ & Other singular fiber \\
\hline
$(lv1)$ & $II*$ & $I_1$ & $ I_1  $\\
$(lv2)$ & $II*$&$I_1$&$ I_1  $\\
$(lv3)$ & $II*$ & $I_1$ & $I_1$ \\
$(lv5)$ & $II*$ & $I_1$ & $I_1$   \\
\hline
$(lv4)$ & $III*$ & $I_2$ &  $ I_1  $\\
\hline
\end{tabular}
\end{center}
\caption{\label{tab5}The Lotka-Volterra case with $3$ singular fibers.}
\end{table}
\begin{table}[h!]
\begin{center}
\begin{tabular}{| c | c | c | }
\hline
Fibration &  $\{t=\infty\}$ &  $t=0$  \\
\hline
$(iii), k=l=1 \mod(3)$ & $IV^*$ & $IV$\\
\ \ \ \ \ $ k=l=2 \mod(3)$ & $IV$ & $IV^*$\\
$(iv), k=l=1 \mod(3)$ &$IV^*$& $IV$\\
\ \ \ \ \ $ k=l=2 \mod(3)$ & $IV$ & $IV^*$\\
\hline
$(v), k=l=1 \mod(4)$ &$III^*$ &$III$ \\
\ \ \ \ \ $ k=l=3 \mod(4)$ &$III$ &$III^*$ \\
$(vi), k=l=1 \mod(4)$ &$III^*$ &$III$ \\
\ \ \ \ \ $ k=l=3 \mod(4)$ &$III$ &$III^*$ \\
\hline
$(vii), l=2 \mod(6)$  & $II^*$ & $II$ \\
\ \ \ \ \ $ l=4 \mod(6)$  & $II$ & $II^*$ \\
$(viii), l=2 \mod(6)$ & $II^*$& $II$ \\
\ \ \ \ \ $ l=4 \mod(6)$  & $II$ & $II^*$\\
\hline
$(ix), k=1 \mod(4)$ & $III^*$ & $III$\\
\ \ \ \ \ $ k=3 \mod(4)$ & $III$ & $III^*$\\
$(x), k=1 \mod(4)$ & $III^*$ & $III$\\
\ \ \ \ \ $ k=3 \mod(4)$ & $III$ & $III^*$\\
\hline
$(xi), k=1,2 \mod(6)$ & $II^*$ & $II$\\
\ \ \ \ \ $k=4,5 \mod(6)$ & $II$ & $II^*$\\
\hline
$(xii), k=1,2 \mod(6)$ & $II^*$ & $II$\\
\ \ \ \ \ $k=4,5 \mod(6)$ & $II$ & $II^*$\\
\hline
\end{tabular}
\end{center}
\caption{\label{tab6} The elliptic Lotka-Volterra case with $2$ singular fibers.}
\end{table}

\newpage

\subsection{The reversible Lotka-Volterra case}
The calculus are similar and are left to the reader. The results are contained in Tables \ref{tab7} and \ref{tab8}.
\begin{table}[h!]
\begin{center}
\begin{tabular}{| c | c | c | c |}
\hline
Fibration & Fiber $\{t=\infty\}$ & Fiber $t=0$ & Other singular fiber \\
\hline
$\rm(rlv1)$ &$IV*$ &$I_3$ &  $I_1$  \\
$\rm(rlv2)$ & $IV*$ &$I_3$ &$I_1$    \\
\hline
$\rm(rlv3)$ &$III*$  & $I_2$& $I_1$   \\
$\rm(rlv4)$ &$III*$ &$I_2$ &   $I_1$ \\
\hline
$\rm(rlv5)$ &$IV$  & $I_1*$& $I_1$   \\
$\rm(rlv6)$ &$IV$ &$I_1*$ &   $I_1$ \\
\hline
\end{tabular}
\end{center}
\caption{\label{tab7} The elliptic reversible Lotka-Volterra case with $3$ singular fibers.}
\end{table}
\begin{table}[h!]
\begin{center}
\begin{tabular}{| c | c | c | }
\hline
Fibration & Fiber $\{t=\infty\}$ & Fiber $t=0$ \\
\hline
$(ix), k=1 \mod(3)$&$IV*$  &$IV$   \\
\ \  \ \ \ $ k=2\mod(3)$&$IV$  &$IV*$   \\
$(x), k=1 \mod(3)$&$IV*$ & $IV$   \\
\ \  \ \ \ $ k=2\mod(3)$&$IV$  &$IV*$  \\
\hline
$(xi), k=1 \mod(4)$&$III*$  &$III$   \\
\ \  \ \ \ $ k=3\mod(4)$&$III$  &$III*$  \\
$(xii), k=1 \mod(4)$&$III*$ &$III$   \\
\ \  \ \ \ $ k=3\mod(4)$&$III$  &$III*$  \\
\hline

\end{tabular}
\end{center}
\caption{\label{tab8} The elliptic reversible Lotka-Volterra case with $2$ singular fibers.}
\end{table}

\vspace{10cm}
\textbf{\large Acknowledgements:}\\
\\
The author would like to thank L. Gavrilov and I. D. Iliev and  for their attention to the paper, many useful comments and stimulating remarks.


\begin{thebibliography}{2}
\bibitem[BPV84] {surfaces} W. Barth , C. Peters, A. Van de Ven: \textsl{Compact complex surfaces}, Erg. Math., Springer-Verlag (1984).
\bibitem[CGLPR]{poly}J. Chavarriga, B. Garci\'a, J. Llibre, J. S P\'erez Del R\'io and J. A. Rodr\'iguez: \textsl{Polynomial first integrals of quadratic vector fields}, Journal of Differential Equations
Volume 230, Issue 2, 15 November 2006, Pages 393-421.
\bibitem[CLLL06]{chen} G. Chen, C. Li, C. Liu, J. Llibre: \textsl{The cyclicity of period annuli of some classes of
reversible quadratic systems}, Discrete Contin. Dyn. Syst. { \bf 16}
(2006), no. 1, 157-177.
\bibitem[G07]{thesis}S. Gautier: \textsl{Feuilletages elliptiques
    quadratiques plans et leurs pertubations}, P.H.D Thesis, 7/12/2007.
\bibitem[GGI]{ggi}S. Gautier, L. Gavrilov, I. D. Iliev:
  \textsl{Pertubations of quadratic centers of genus one}, http://arxiv.org/abs/0705.1609.
\bibitem[G01]{hamiltonian}L. Gavrilov: \textsl{The infinitesimal 16th Hilbert problem in the quadratic case}, Invent. math 143, 449-497 (2001).
\bibitem[GH78]{algebraic} P. Griffiths and J. Harris: \textsl{Principles of algebraic geometry}, John Wiley and Sons, New York (1978).
\bibitem[I98]{centers}I. D. Iliev: \textsl{Pertubations of quadratic centers}, Bull. Sci. math. 22 (1998). 107-161.
\bibitem[ILY05]{iliev} I. D. Iliev, C.Li, J.Yu: \textsl{Bifurcations of limit cycles
    from quadratic non-Hamiltonian systems with two centres and two unbounded
    heteroclinic loops}, Nonlinearity, 18 (2005), no. 1, 305--330.
\bibitem[I02]{Ilyashenko} Y. Ilyashenko: \textsl{Centennial history of Hilbert's 16th problem}, Amer. Math. So., 39 (2002), no 3, 301-354.
\bibitem[J79]{pfaff}JC. Jouanolou: \textsl{Equations de Pfaff alg\'ebriques}, Lec. Notes in Math., 708 (1979).
\bibitem[Ka75]{K75} K. Kas: \textsl{Weirstrass normal forms and invariants of elliptic surfaces}, Trans. of the A.M.S, 225, (1977).
\bibitem[Ko60]{singular}K. Kodaira: \textsl{On compact analytic surfaces}, I, Annals of Math., 71 (1960),111-152.
\bibitem[Ko63]{singularbis}K. Kodaira: \textsl{On compact analytic
    surfaces}, II, Annals of Math., 77 (1963),563-626.
\bibitem[M89]{mir}R. Miranda, \textsl{The basic theory of elliptic
    surfaces}, Dottorato di Ricerca in Matematica. [Doctorate in
  Mathematical Research], ETS Editrice, Pisa, 1989.
\bibitem[P90]{Petrov}G. S. Petrov: \textsl{Nonoscillation of elliptic integrals}, Funct. Anal. Appl., 3 (1990), 45-50
\bibitem[YL02]{yu} J. Yu, C. Li: \textsl{Bifurcation of a class of planar
    non-Hamiltonian integrable systems with one center and one homoclinic loop},
    J. Math. Anal. Appl., 269 (2002), no. 1, 227--243.
\end{thebibliography}
\end{document}